\documentclass[a4paper,11pt]{article}
\usepackage[cp1251]{inputenc}
\usepackage[english, russian]{babel}
\usepackage{amsmath,layout}
\usepackage{amsfonts}
\usepackage{amssymb}
\usepackage{amsthm}
\usepackage{graphicx}
\usepackage{epstopdf}
\usepackage{wrapfig}
\usepackage{verbatim}
\usepackage{float}
\usepackage{indentfirst}
\usepackage{graphicx}

\usepackage{amsfonts,amssymb,mathrsfs,amscd}

\newcommand{\sgrad}[0]{\mathrm{sgrad \,}}

\theoremstyle{plain}
\newtheorem{theorem}{Теорема}

\newtheorem{corollary}{Следствие}[section]

\theoremstyle{definition}
\newtheorem{definition}{Определение}[section]

\newtheorem{example}{Пример}[section]

\newtheorem{coment}{Замечание}[section]

\def\Re{\operatorname{Re}}
\def\Im{\operatorname{Im}}
\def\const{\operatorname{const}}
\def\RR{\mathbb R}
\def\CC{\mathbb C}
\def\NN{\mathbb N}

\def\ZZ{\mathbb Z}

\def\sgrad{{\rm sgrad\,}}

\def\conv{{\rm conv}}

\def\Pr{{\rm Pr}}
\def\SSigma{{\rm T}}

\def\eps{\varepsilon}
\def\Spec{{\rm Spec}}

\def\ZZ{\mathbb Z}
\def\NN{\mathbb N}
\def\RR{\mathbb R}
\def\CC{\mathbb C}

\def\NP{{\rm P}}

\def\Re{\mathop{\mathrm{Re}}\nolimits}
\def\Im{\mathop{\mathrm{Im}}\nolimits}
\def\sgrad{{\rm sgrad\,}}

\def\conv{{\rm conv}}

\def\Pr{{\rm Pr}}

\def\Sym{{\rm Sym}}

\def\<{\langle}
\def\>{\rangle}
\def\eps{\varepsilon}
\def\d{\partial}

\title{Integrable Hamiltonian systems with incomplete flows \\ and Newton's polygons}
\author{Elena A.\ Kudryavtseva, Timur A.\ Lepsky}
\date{}

\begin{document}

\maketitle

{\small {\bf Abstract.} We study the Hamiltonian vector field $v=(-\partial f/\partial w,\partial f/\partial z)$ on $\mathbb C^2$, where $f=f(z,w)$ is a polynomial in two complex variables, which is non-degenerate with respect to its Newton's polygon. We introduce coordinates in four-dimensional neighbourhoods of the ``points at infinity'', in which the function $f(z,w)$ and the 2-form $dz\wedge dw$ have a canonical form. A compactification of a four-dimensional neighbourhood of the non-singular level set $T_0=f^{-1}(0)$ of $f$ is constructed. The singularity types of the vector field
$v|_{T_0}$ at the ``points at infinity'' in terms of Newton's polygon are determined.


{\it MSC-class: } 37J05, 37J35
}

 \begin{center} \bf {\LARGE
Интегрируемые гамильтоновы системы с неполными потоками и многоугольники Ньютона} \medskip\\
Елена А.\ Кудрявцева, Тимур А.\ Лепский\footnote{Работа выполнена при поддержке
Российского фонда фундаментальных исследований (гранты №№
10--01--00748-а, 08--01--91300-ИНД a), Программы поддержки ведущих
научных школ РФ (грант № НШ-3224.2010.1), Программы ``Развитие
научного потенциала высшей школы'' (грант № 2.1.1.3704) и Программы
ФЦП <<Научные и научно-педагогические кадры инновационной России>>
(грант № 02.740.11.5213).}
 \end{center}

\begin{abstract}
Изучается гамильтоново векторное поле $v=(-\partial f/\partial w,\partial f/\partial z)$ в $\mathbb C^2$, где $f=f(z,w)$ -- многочлен двух комплексных переменных, невырожденный относительно своего многоугольника Ньютона. Введены координаты в четырехмерных окрестностях ``бесконечно
удаленных точек'', в которых функция $f(z,w)$ и 2-форма $dz\wedge
dw$ имеют канонический вид. Построена компактификация четырехмерной окрестности
неособого множества уровня $T_0=f^{-1}(0)$ функции $f$. Вычислены типы особенностей
векторного поля $v|_{T_0}$ в ``бесконечно удаленных точках'' в терминах многоугольника Ньютона.

\medskip
{\it Ключевые слова: } интегрируемая гамильтонова система, неполные гамильтоновы потоки, многоугольник Ньютона.

{\it УДК } 517.938.5, 514.756.4
\end{abstract}


\section{Введение}

В работе рассматривается многочлен $f:\CC^2\to\CC$ двух комплексных переменных, отличный от константы.
Согласно наблюдению А.И.~Шафаревича, для такой функции $f$
гамильтонова система $(\CC^2,\Re(dz\wedge dw),\Re f(z,w))$ обладает дополнительным первым интегралом $\Im f$ и
равносильна $\CC$-гамильтоновой системе $(\CC^2,dz\wedge dw,f(z,w))$
(см.\ определение~\ref {newton:def:0.4} или~\cite{BaCu,leps}). Данные
гамильтоновы системы часто являются
{\it гамильтоновыми системами с неполными потоками}, т.е.\ пара
коммутирующих векторных полей $\sgrad\Re f=(-\frac{\d f}{\d
w},\frac{\d f}{\d z})$ и $\sgrad\Im f$ (или $\sgrad_\CC f=(-\frac{\d
f}{\d w},\frac{\d f}{\d z})$ и $i\,\sgrad_\CC f$) обладает неполными
потоками на любом слое
 $$
\SSigma_{\xi}:=f^{-1}(\xi)=\{(z,w)\in\CC^2\mid f(z,w)=\xi\}, \quad
\xi\in\CC.
 $$
Отдельно отметим, что любой слой системы не является компактным,
причем при некоторых естественных ограничениях на функцию $f$ почти
все слои являются неособыми и связными (см.\ теоремы~\ref
{newton:lem:dis} и~\ref {newton:the:1.2}). Тем самым, классическая
теорема Лиувилля неприменима для таких интегрируемых гамильтоновых
систем. Возникает задача (поставленная А.Т.~Фоменко)
об обобщении классической теоремы Лиувилля на случай интегрируемых
гамильтоновых систем с неполными потоками, а именно: описание
топологии неособого слоя, описание лагранжева слоения в окрестности неособого
слоя, построение аналога координат действие-угол. Мы рассматриваем
эту задачу для указанного выше класса интегрируемых гамильтоновых систем с неполными потоками, который был введен Х.\ Флашкой~\cite{Flaschka} и предложен авторам А.Т.~Фоменко и А.И.~Шафаревичем. В случае гиперэллипитических многочленов $f(z,w)=z^2+P(w)$, где $P(w)$ -- многочлен с простыми вещественными корнями, задача решена авторами в~\cite {leps,KLmsb,KLtsvta}.

В настоящей работе определено два способа пополнения исходной системы и ее
ограничения на слой. Первый способ заключается в компактификации
многообразия $\bigcup\limits_{|\xi-\xi_0|\le\eps}\SSigma_\xi$ с
краем путем добавления к каждому слою $\SSigma_\xi$ некоторого числа
``бесконечно удаленных точек'' и описания их четырехмерных
окрестностей, см.\ теорему~\ref{newton:intrth2} и следствие~\ref
{newton:cor2} (этот способ пополнения аналогичен теореме
А.Г.~Хованского~\cite[\S2.2, теорема и замечание 1]{hov.1} о
разрешении особенностей при помощи подходящего торического
многообразия). При таком пополнении компактификация
$\widetilde\SSigma_\xi$ любого неособого слоя $\SSigma_\xi$ данного
семейства является неособой компактной связной двумерной вещественной
поверхностью.

Второй способ заключается в пополнении индивидуального слоя
$\SSigma_\xi$, $\xi\in\CC$, системы относительно метрики на слое
$\SSigma_\xi$, естественным образом связанной с векторным полем
$\sgrad_\CC f|_{\SSigma_\xi}$, см.\ определение~\ref
{newton:def:0.6} и следствие~\ref{newton:cor:Newton}. Отметим, что
результатом пополнения слоя $\SSigma_\xi$ вторым способом является
связная двумерная вещественная поверхность $\overline\SSigma_\xi$,
которая не обязательно компактна, причем
ее компактность равносильна неполноте потоков пары векторных полей
$\sgrad_\CC f|_{\SSigma_\xi},i\,\sgrad_\CC f|_{\SSigma_\xi}$, а
именно: в случае неполноты потоков пополнение слоя компактно
(например, для ненулевых слоев функции $f(z,w)=z^2+w^n$ при $n\ge3$),
а в случае полноты потоков пополнение слоя совпадает с самим слоем и
потому некомпактно (например, для ненулевых слоев функции
$f(z,w)=z^2+w^n$ при $n=1,2$).
Более точно: результатом пополнения слоя $\SSigma_\xi$ вторым
способом является связная двумерная вещественная поверхность
$\overline\SSigma_\xi$ с плоской римановой метрикой и конечным
числом конических особенностей, с углом в каждой кратным полному и
выражающимся в терминах многоугольника Ньютона, причем поверхность
$\overline\SSigma_\xi$ либо компактна и гомеоморфна сфере с
положительным числом ручек, либо некомпактна и изометрична
евклидовой плоскости или плоскому цилиндру (см.\
следствие~\ref{newton:cor:Newton}). Поток в силу данной
гамильтоновой системы на такой компактной комплексной кривой
$\overline\SSigma_\xi$ определяется голоморфной 1-формой
$\Delta_\xi$ на $\overline\SSigma_\xi$ (см.\ определение~\ref
{newton:def:0.6}(А)); динамические свойства таких потоков изучались,
например, в работах М.Л.~Концевича-А.В.~Зорича, С.П.~Новикова,
А.И.~Буфетова и др. (см.\ работы~\cite
{konts},~\cite{nov},~\cite{buf} и ссылки в них).

Другим результатом работы является описание в терминах многоугольника
Ньютона типов особенностей векторного поля $\sgrad_\CC
f|_{\SSigma_\xi}$ в ``бесконечно удаленных точках'', см.\
следствия~\ref {newton:cor2} и~\ref{newton:cor:Newton}.

В \S\ref{sec:0} и \S\ref{sec:0.1} содержится обзор известных
результатов
(см.~\cite{hov.1},~\cite{hov.2},~\cite{hart},~\cite{forster}) о
конечности множества критических значений и о топологии нулевого слоя
$\SSigma_0=f^{-1}(0)$ для многочлена $f$.
Появление данного обзора мотивировано тем, что во многих работах, где
изучаются топологические свойства алгебраических множеств (т.е.\
совместных множеств уровня многочленов), результаты сформулированы
для многочлена (или системы многочленов) ``общего положения'',
без явного перечисления условий (аналогичных условиям типа
невырожденности, перечисленным в определениях~\ref {def:nondegen},
\ref {def:nonsingular}, \ref {def:irreducible}), накладываемых на
многочлен. Это объясняется тем, что, в частности, многочлены, не
удовлетворяющие данным условиям, образуют множество меры нуль в
пространстве всех многочленов. Однако такой подход лишает возможности
использования подобных теорем в приложениях, поскольку не дает ответа
на вопрос о применимости таких результатов к конкретным многочленам.
В настоящем обзоре сформулированы достаточные условия, накладываемые
на многочлены, для возможности применения таких теорем на практике.

Применение результатов работы к конкретным многочленам иллюстрируется
в~\S\ref {sec:5}, а также в примере~\ref {newton:ex:disconnected}.

Авторы выражают глубокую признательность А.Т.~Фоменко за постановку задачи,
А.Т.~Фоменко и А.И.~Шафаревичу за предложенную для исследования систему,
А.Б.~Жеглову за многочисленные обсуждения.

\section{Основные понятия и утверждения} \label {sec:0}

\begin{figure}[h]
\begin{center}
\includegraphics[scale=0.6]{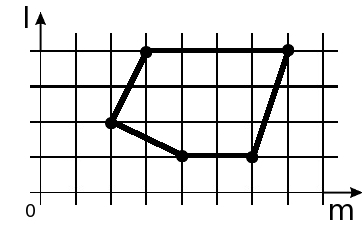}
\caption{Пример многоугольника Ньютона} \label{n:fig:ex_g}
\end{center}
\end{figure}

\begin{definition} [{\rm \cite[\S2.1 и \S3.4]{hov.1}}] \label{newton:def:0.1}
{\em Многоугольником Ньютона} $\NP_f\subset\RR^2$ многочлена
$f(z,w)=\sum\limits_{l,m\ge0}a_{l,m}z^lw^m$ называется выпуклая
оболочка множества точек $(l,m)\in\ZZ^2$ таких, что $a_{l,m}\ne0$.
{\em Размерностью} $\dim\NP_f$ многоугольника Ньютона $\NP_f$
называется размерность минимального аффинного подпространства в
$\RR^2$, содержащего многоугольник Ньютона $\NP_f$.
{\em Многоугольником} $\NP_f^{\eta}\subset\RR^2$, отвечающим
многоугольнику Ньютона $\NP_f$ и ковектору $\eta\in\RR^{2*}$,
назовем грань (размерности $0,1$ или $2$) многоугольника Ньютона
$\NP_f$, на которой достигает максимума функция $\NP_f\to\RR$,
$x\mapsto\langle\eta, x\rangle$, $x\in\NP_f$, где через
$\langle\eta, x\rangle\in\RR$ обозначено значение ковектора $\eta$
на векторе $x\in\RR^2$ (в частности, сам многоугольник Ньютона
$\NP_f^{0}=\NP_f$ отвечает нулевому ковектору $\eta=0$, а каждая
сторона многоугольника $\NP_f$ отвечает своему ``вектору внешней
нормали''). По многочлену
$f(z,w)=\sum\limits_{l,m\ge0}a_{l,m}z^lw^m$ и ковектору $\eta$
определим {\em усеченный многочлен}
$f^{\eta}(z,w)=\sum\limits_{(l,m)\in\NP_f^{\eta}}a_{l,m}z^lw^m$, см.\
рис.~\ref{n:fig:ex_g}.
\end{definition}

Ниже (определения~\ref {def:nondegen}, \ref {def:nonsingular} и~\ref
{def:irreducible}) введены три понятия невырожденности для многочлена
$f(z,w)-\xi$, а именно определены {\it невырожденность многочлена
относительно своего многоугольника Ньютона}, {\it неособость} слоя
$f^{-1}(\xi)$ {\em для функции} $f$ и {\it неприводимость} многочлена
$f(z,w)-\xi$. Любое из этих условий выполнено для многочленов
``общего положения''. При выполнении данных условий удается описать
топологические свойства слоя $f^{-1}(\xi)$ (см.\ теоремы~\ref
{newton:the:1.2}, \ref {newton:the:1.1}, \ref {newton:the:1.3}, \ref
{newton:the:1.4} и следствие~\ref {newton:cor:1.3}, а также
теорему~\ref {newton:intrth2} и следствия~\ref {newton:cor2} и~\ref
{newton:cor:Newton} ниже). Как будет показано ниже, данные свойства
являются независимыми (см.\ пример~\ref {newton:ex:types}).

\begin{definition} [{\rm \cite[\S2.1, определение]{hov.1}}] \label{def:nondegen}
Многочлен $f(z,w)$ называется {\em невырожденным относительно своего
многоугольника Ньютона} $\NP_f$, если для любого ковектора
$\eta\in\RR^{2*}$ выполнено следующее условие: для любого решения
$(z,w)$ уравнения $f^{\eta}(z,w)=0$, лежащего в
$(\CC\setminus\{0\})^2$, дифференциал $df^{\eta}(z,w)\ne0$.
\end{definition}

\begin{definition}\label{newton:def:0.5} \label{def:nonsingular}
{\it Слоем} $\SSigma_{\xi}\subset\CC^2$ многочлена
$f:\CC^2(z,w)\to\CC$ назовем множество
$\SSigma_{\xi}=\{(z,w)\in\CC^2|f(z,w)=\xi\}$.
Слой $\SSigma_{\xi}\subset\CC^2$ называется {\it неособым} для
функции $f$, если $df(z,w)\ne0$ для любых $(z,w)\in\SSigma_{\xi}$.
\end{definition}

\begin{definition} \label{def:irreducible}
Многочлен $f(z,w)$ называется {\it неприводимым}, если не существует
его разложения на множители $f(z,w)=p_1(z,w)p_2(z,w)$, где $p_1(z,w)$
и $p_2(z,w)$ --- многочлены, отличные от константы.
\end{definition}

В дальнейшем в основном будут рассматриваться только неособые слои
невырожденных неприводимых многочленов. Естественность этого
предположения показывает следующее замечание~\ref{newton:com:types}.

\begin{coment} \label{newton:com:types}
Для заданного многоугольника $\NP\ni(0,0)$ размерности $\dim\NP=2$ и для почти всех многочленов $f$, таких что $\NP_{f}=\NP$, выполнены следующие свойства:

1) нулевой слой $\SSigma_0=\{(z,w)\in\CC^2\mid f(z,w)=0\}$ является
неособым (см.\ теорему~\ref {newton:lem:dis} или~\ref {newton:the:1.3});

2) многочлен $f(z,w)$ является невырожденным относительно своего
многоугольника Ньютона $\NP_f$ (см.~\cite [\S2.2, теорема]{hov.1});

3) многочлен $f(z,w)$ является неприводимым (см.\ теорему~\ref {newton:the:1.1}).
\end{coment}

\begin{example}\label{newton:ex:types}
В таблице ниже для каждого из выписанных многочленов указано,
выполнены ли для него свойства неособости нулевого слоя,
невырожденности относительно своего многоугольника Ньютона и
неприводимости.
\begin{center}
\begin{tabular}{lccc}
Многочлен $f(z,w)$ &Неособость $\SSigma_0$&Невырожденность $f$&Неприводимость $f$\\
$z$&$+$&$+$&$+$\\
$z^2+w^3$&$-$&$+$&$+$\\
$z^2+w^2+2zw+w$&$+$&$-$&$+$\\
$z^3+(w+1)^2$&$-$&$-$&$+$\\
$(z+1)(z+2)$&$+$&$+$&$-$\\
$zw$&$-$&$+$&$-$\\
$(z^2+w^2+1)(z^2+w^2+2)$&$+$&$-$&$-$\\
$(z+1)^3$&$-$&$-$&$-$
\end{tabular}
\end{center}
\end{example}

\begin{coment}
Как показано в примере~\ref{newton:ex:types}, условия неособости слоя
$\SSigma_\xi$ многочлена $f$, невырожденности многочлена $f-\xi$
относительно своего многоугольника Ньютона $\NP_{f-\xi}$,
неприводимости многочлена $f-\xi$ являются независимыми, где
$\xi\in\CC$.
\end{coment}

\begin{theorem}[Конечность множества особых значений {\rm\cite{hart}, \cite{mats}}]\label{newton:lem:dis}
Пусть $f:\CC^2\to\CC$ --- комплексный многочлен двух комплексных
переменных, отличный от константы. Тогда множество
$\Sigma_f\subset\CC$ особых значений конечно, т.е.\ имеет вид
$\Sigma_f=\{\xi_i\}_{i=1}^N$, где $\xi_i\in\CC,\,i=1,\dots,N$.
\end{theorem}

\begin{proof}
Рассмотим в $\CC^3$ гладкое подмногообразие
$X=\{(z,w,\xi)\in\CC^3\mid f(z,w)=\xi\}\subset\CC^3$ и регулярное
отображение $\Pr_\xi:\CC^3\to\CC$, $(z,w,\xi)\mapsto \xi$.
Регулярность отображения $\Pr_\xi$ следует из того, что, в частности,
$\frac{\d\Pr_\xi(z,w,\xi)}{\d \xi}=1$. Образом отображения
$F:=\Pr_\xi|_{X}$ является $F(X)=\CC$, поскольку многочлен $f(z,w)$
отличен от константы. По усиленной теореме Бертини (для схем) отсюда
следует (см.~\cite[Гл.III, следствие 10.7]{hart}), что существует
открытое по Зарисскому непустое множество $U\subset\CC$ такое, что
$F|_{F^{-1}(U)}:F^{-1}(U)\to U$ является гладким морфизмом
соответствующих схем
в смысле~\cite[Гл.III, \S10, определение]{hart}.
По теореме~\cite[Гл.III, теорема 10.2]{hart} схема
$X_\xi:=\Spec(\CC[z,w]/(f-\xi))$ является регулярной для любого
$\xi\in U$, т.е., в частности,
она регулярна в следующем смысле: для любой точки $(z,w)\in\CC^2$,
такой что $(z,w,\xi)\in X$, локальное кольцо данной схемы в точке
$(z,w,\xi)$ регулярно в смысле~\cite[Гл.I, \S5, определение]{hart}.
Покажем, что отсюда
следует, что $(\frac{\d f}{\d z}(z,w),\frac{\d f}{\d
w}(z,w))\ne(0,0)$, т.е.\ $\xi$ является неособым значением функции
$f$. Действительно, из теоремы (см.~\cite[теорема 36, с.121]{mats}
или~\cite[Гл.1, упражнение 5.13]{hart}) о том, что локальное
регулярное кольцо не имеет делителей нуля, следует, что $f(z,w)-\xi$
--- это произведение неприводимых многочленов $P_i(z,w)$ без общих
нулей, поэтому неравенство $(\frac{\d f}{\d z}(z,w),\frac{\d f}{\d
w}(z,w))\ne(0,0)$ равносильно системе аналогичных неравенств для
каждого неприводимого сомножителя $P_i(z,w)$. А для неприводимого
многочлена $P(z,w)$ требуемое неравенство доказано в~\cite[Гл.I,
теорема 5.1]{hart}.
Всякое открытое по Зарисскому непустое подмножество $U\subset\CC$
имеет вид $U=\CC\setminus\{\xi_i\}_{i=1}^{N}$, откуда образ
$\Sigma_f$ множества особых точек содержится в $\{\xi_i\}_{i=1}^N$,
т.е.\ конечен.
\end{proof}

\section{Обзор известных результатов по топологии слоев} \label {sec:0.1}

\subsection{Достаточные условия связности слоя} \label {subsec:0.1}

Следующие две теоремы~\ref{newton:the:1.2} и~\ref {newton:the:1.1}
устанавливают связность нулевого слоя многочлена, являющегося либо
неприводимым, либо невырожденным относительно своего многоугольника
Ньютона, имеющего размерность 2.

\begin{theorem} [Связность слоя неприводимого многочлена {\rm \cite[Гл.I]{hart}}] \label{newton:the:1.2}
Пусть многочлен $f=f(z,w)$ неприводим. Тогда нулевой слой $\SSigma_0=f^{-1}(0)$ связен.
\end{theorem}

\begin{proof} Согласно~\cite[Гл.I, следствие 1.4]{hart} слой $\SSigma_0$ неприводим тогда и только тогда, когда идеал, порожденный многочленом
$f$, является простым. Это, в свою очередь эквивалентно тому, что
многочлен $f$ является неприводимым. Неприводимость слоя
$\SSigma_0$ означает (см.~\cite[Гл.I, определение в \S1, стр.
18]{hart}), что не существует $Y_1,Y_2\subset\SSigma_0$ ---
собственных замкнутых (в смысле топологии Зарисского) в $\SSigma_0$
подмножеств, таких что $\SSigma_0=Y_1\bigcup Y_2$, что влечет
связность $\SSigma_0$ в смысле топологии Зарисского.
Замыкание $\overline{\SSigma}_0$ слоя $\SSigma_0$ в проективной плоскости
также является неприводимым (см.~\cite[Гл.1, \S1, упражнение]{hart}), причем $\overline{\SSigma}_0\setminus\SSigma_0$ ---
конечное множество точек. Отсюда проективная кривая $\overline{\SSigma}_0$ связна в топологии Зарисского. По теореме Серра (см.~\cite[Добавление C]{hart}) кривая $\overline{\SSigma}_0$ линейно связна, откуда
слой $\SSigma_0$ линейно связен.
\end{proof}

\begin{theorem} [Связность слоя невырожденного многочлена {\rm \cite[\S2.1]{hov.2}}] \label{newton:the:1.1}
Пусть многочлен $f=f(z,w)$ невырожден относительно своего
многоугольника Ньютона $\NP_f$, и $\dim\NP_f=2$. Тогда подмножество
$\hat\SSigma_0:=\SSigma_0\setminus((\CC\times\{0\})\cup(\{0\}\times\CC))$
нулевого слоя $\SSigma_0=f^{-1}(0)$ связно.
\end{theorem}

\begin{proof} Пусть $\widetilde X$ --- замыкание многообразия
$X:=\hat\SSigma_0$ в достаточно полной проективной торической
компактификации $M\supset\CC^2$. Согласно~\cite[\S2.1,
теорема]{hov.2}, $\widetilde X$ связно. Так как $\widetilde X,X$
бирационально эквивалентны и имеют комплексную размерность 1,
множество $\widetilde X\setminus X$ конечно. Поэтому из связности
$\widetilde X$ получаем связность $X$.
\end{proof}

\subsection{Топология слоя невырожденного многочлена} \label {subsec:0.2}

\begin{definition} \label{def:newton:arithm}
Пусть $M$ --- компактное аналитическое проективное многообразие
размерности $\dim_\CC M=n$. {\it Арифметическим родом} $p_a(M)$
многообразия $M$ называется альтернированная сумма
$$
p_a(M):=\sum_{k=0}^n(-1)^k\dim_\CC(\Omega^k(M)),
$$
где через $\Omega^k(M)$ обозначено пространство голоморфных
дифференциальных $k$-форм на $M$ (см.~\cite[\S1.1]{hov.2}). Для
некомпактного аналитического многообразия $M'$ арифметический род
определяется формулой $p_a(M'):=p_a(M)$, где $M$ --- любое компактное
аналитическое проективное многообразие, бирационально эквивалентное
$M'$. \end{definition}

Известно, что для бирационально эквивалентных компактных
алгебраических многообразий $M_1,M_2$ выполнено
$\dim_\CC(\Omega^k(M_1))=\dim_\CC(\Omega^k(M_2))$ и
$p_a(M_1)=p_a(M_2)$ (см.~\cite[\S1.1]{hov.2}). Поэтому
определение~\ref{def:newton:arithm} корректно для некомпактных
аналитических многообразий. Если $M_1$ и $M_2$ бирационально
эквивалентны и $\dim_\CC M_1=1$, то имеются конечные подмножества
$N_1\subset M_1$ и $N_2\subset M_2$, такие что $M_1\setminus N_1$ и
$M_2\setminus N_2$ комплексно диффеоморфны.

\begin{theorem} [Неособость и арифметический род слоя {\rm \cite[\S1.1, теорема 1 в \S1.3, теорема в \S2.1, теорема в \S4.1]{hov.2}}] \label{newton:the:1.3}
Слой $\SSigma_0=f^{-1}(0)$, определенный невырожденным (относительно
своего многоугольника Ньютона $\NP_f$) многочленом $f(z,w)\ne\const$
с ненулевым свободным членом, является неособым для функции $f$.
Его арифметический род $p_a(\SSigma_0)$ вычисляется по формуле
$p_a(\SSigma_0)
= 1-(-1)^{\dim\NP_f}B^+(\NP_f)$, где $B^+(\NP_f)$ --- количество
целочисленных точек, лежащих строго внутри многоугольника Ньютона
$\NP_f$ (в топологии минимального линейного пространства, содержащего
$\NP_f$).
То есть, $p_a(\SSigma_0)=1-B^+(\NP_f)$ при $\dim\NP_f=2$, и
$p_a(\SSigma_0)=1+B^+(\NP_f)$ при $\dim\NP_f=1$.
\end{theorem}

\begin{example} \label {newton:ex:disconnected}
Для $f(z,w)=z^n-1$ слой
$\SSigma_0=f^{-1}(0)\approx\CC\times\{1,\dots,n\}$ имеет
арифметический род $p_a(\SSigma_0)=1+B^+(\NP_f)=n$ в силу
теоремы~\ref {newton:the:1.3}. Соответствующая компактная
аналитическая проективная кривая $M\supset\SSigma_0$, биголоморфно
эквивалентная слою $\SSigma_0$, является несвязным объединением $n$
экземпляров сферы Римана: $M=\overline\CC\times\{1,\dots,n\}$.
\end{example}

Согласно~\cite[\S19.14]{forster}, если $X$ --- компактное связное
1-мерное комплексное многообразие (т.е.\ {\it риманова поверхность}),
гомеоморфное сфере с $n_g$ ручками, то размерность пространства
голоморфных 1-форм на нем равна $n_g$. Отсюда получаем следующее
следствие теорем~\ref {newton:the:1.1} и~\ref {newton:the:1.3}.

\begin{corollary} [Количество ручек у слоя] \label{newton:cor:1.3}
Пусть многочлен $f(z,w)-\xi$ невырожден относительно своего многоугольника Ньютона $\NP_{f-\xi}$,
причем $\dim\NP_{f-\xi}=2$ и $f(0,0)\ne\xi$.
Тогда слой $\SSigma_\xi$ является неособым для функции $f$ и
гомеоморфен сфере с $n_g=B^+(\NP_{f-\xi})$ ручками и конечным числом
проколов, где $B^+(\NP_{f-\xi})$ как в теореме~\ref{newton:the:1.3}.
\end{corollary}

\begin{definition}\label{newton:def:0.4}
Векторным полем {\it косой градиент} $\sgrad_\CC f\in
\rm{Vect}(\CC^2)$ голоморфной функции $f:\CC^2\to\CC$ относительно
голоморфной 2-формы $\omega_\CC=dz\wedge dw$ называется векторное поле $\sgrad_\CC f =
(-\frac{\partial f(z,w)}{\partial w}, \frac{\partial f(z,w)}{\partial
z})$, заданное в координатах $(z,w)$.
\end{definition}

\begin{definition}\label{newton:def:0.6}
(А) {\em Римановой метрикой пополнения} $g_{\xi}$ неособого слоя
$\SSigma_{\xi}$ для функции $f$ назовем риманову метрику
$g_{\xi}=\Sym(\Delta_{\xi}\otimes\overline{\Delta}_{\xi})$, где {\em
голоморфная 1-форма} $\Delta_{\xi}$ определена на слое
$\SSigma_{\xi}$ соотношением $\Delta_{\xi}(\sgrad_\CC
f|_{\SSigma_{\xi}}) = 1$. Отметим, что риманова метрика $g_{\xi}$
является плоской, и интегральные траектории векторных полей
$\sgrad_\CC f|_{\SSigma_{\xi}}$ и $i\,\sgrad_\CC f|_{\SSigma_{\xi}}$
являются ее геодезическими.

(Б) На слое $\SSigma_{\xi}$ определена {\em функция расстояния}
$\rho_{\xi}:\SSigma_{\xi}\times\SSigma_{\xi}\to\mathbb{R}$, где для
любых $x,y\in\SSigma_{\xi}$, $\rho_{\xi}(x,y)$ --- нижняя грань длин
всех кривых, лежащих в $\SSigma_{\xi}$ и соединяющих точки $x,y$,
расстояние в смысле римановой метрики пополнения $g_{\xi}$.
\end{definition}

\begin{theorem} [Количество ручек и голоморфные 1-формы на слое {\rm \cite[утверждение и пример в \S2.2]{hov.2}}]  \label{newton:the:1.4}
Пусть многочлен $f(z,w)-\xi$ невырожден относительно своего
многоугольника Ньютона $\NP_{f-\xi}$, причем $\dim\NP_{f-\xi}=2$.
Тогда подмножество
$\hat\SSigma_\xi:=\SSigma_\xi\setminus((\CC\times\{0\})\cup(\{0\}\times\CC))$
слоя $\SSigma_\xi=f^{-1}(\xi)$
гомеоморфно сфере с $n_g=B^+(\NP_{f-\xi})$ ручками и конечным числом
проколов, где $B^+(\NP_{f-\xi})$ как в теореме~\ref{newton:the:1.3}.
Более того, 1-формы $\Delta_\xi z^lw^m$ на $\hat\SSigma_\xi$, где $(l,m)\in\ZZ^2$ ---
внутренние целочисленные точки многоугольника Ньютона $\NP_{f-\xi}$,
образуют базис пространства голоморфных 1-форм на некоторой
компактной связной аналитической проективной кривой $\widetilde
X\supset\hat\SSigma_\xi$, биголоморфно эквивалентной многообразию
$\hat\SSigma_\xi$.
\end{theorem}

\section{Поведение гамильтонова поля в бесконечно удаленных точках на пополненном слое} \label {sec:1}

\begin{definition}\label{newton:def:0.7}
(А) Скажем, что мероморфное векторное поле $v$, определенное на
некоторой комплексной кривой, имеет {\em полюс порядка $k\ge0$ в
точке $x$}, если в некоторой окрестности $U$ точки $x$ выполнено
соотношение: $v = h(u)u^{-k}\frac{d}{du}$, где $u: U\to \CC$ ---
локальная координата в окрестности точки $x$, $h(u)$ --- некоторая
голоморфная функция на $u(U)$, такие что $u(x)=0$ и $h(0)\ne0$.

(Б) Скажем, что голоморфная 1-форма $\Delta$, определенная на
некоторой комплексной кривой, имеет {\em ноль порядка $k\ge0$ в точке
$x$}, если в некоторой окрестности $U$ точки $x$ выполнено
соотношение: $\Delta=h(u)u^kdu$, где $u: U\to \CC$ --- локальная
координата в окрестности точки $x$, $h(u)$ --- некоторая голоморфная
функция на $u(U)$, такие что $u(x)=0$ и $h(0)\ne0$.
\end{definition}

\begin{figure}[h]
\begin{center}
\includegraphics[scale=0.4]{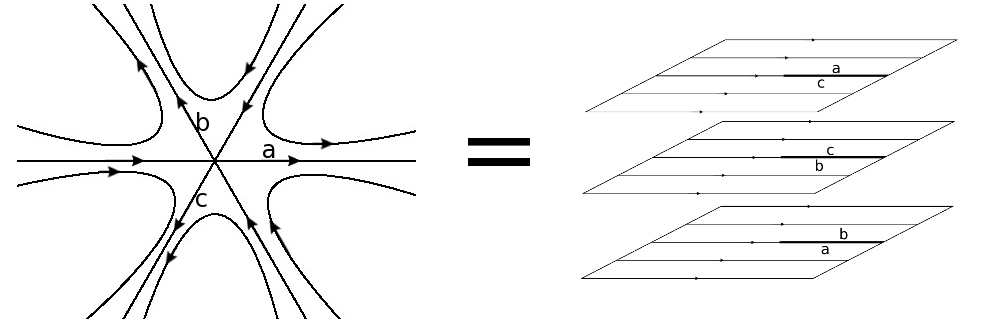}
\caption{Полюс порядка 2} \label{newton:fig:bif_re}
\end{center}
\end{figure}

Из определения~\ref{newton:def:0.7} легко следует, что если
голоморфная 1-форма $\Delta$ на комплексной кривой имеет ноль порядка
$k$ в точке $x$, то в некоторой окрестности $U$ точки $x$ выполнено
соотношение: $\Delta=u^kdu$ для некоторой локальной координаты $u:
U\to \CC$ в окрестности точки $x$. Полюс порядка $k=0$ является
устранимой особенностью векторного поля. Интегральные траектории
векторного поля $v$, имеющего полюс порядка 2, изображены на
рис.~\ref{newton:fig:bif_re}.

Пусть многоугольник Ньютона $P_{f-\xi_0}$ удовлетворяет следующему
условию:
\\
(i) многоугольник Ньютона $P=P_{f-\xi_0}$ содержит вместе с каждой
своей точкой $(u,v)\in P$ прямоугольник
$\conv\{(0,0),(u,0),(0,v),(u,v)\}=P_{z^uw^v+z^u+w^v+1}$.

Условие (i) эквивалентно тому, что $P_{z^l+w^m+1}\subseteq
P_{f-\xi_0}\subseteq P_{z^lw^m+z^l+w^m+1}$ для некоторых
неотрицательных $l,m\in\ZZ$.

Обозначим $D^2_{z_0,\eps}:=\{z\in\CC\mid|z-z_0|<\eps\}$, открытый
двумерный диск.

\begin{theorem}[Нормализация невырожденного многочлена и 2-формы $dz\wedge dw$ в ``бесконечно удаленных точках'' слоев] \label{newton:intrth2}
Пусть $f(z,w)-\xi_0$ --- невырожденный многочлен относительно своего
многоугольника Ньютона $\NP_{f-\xi_0}$, причем многоугольник Ньютона $\NP_{f-\xi_0}$ удовлетворяет условию {\rm (i)} выше, и $\dim\NP_{f-\xi_0}=2$ {\rm (см.\ определение~\ref
{newton:def:0.1})}. Тогда существуют $\eps>0$ и $R>0$, такие что

1) для любой стороны $\Gamma_l$ многоугольника Ньютона, не лежащей на
координатных осях,
существуют ровно $n_{\Gamma_l}$ голоморфных вложений
$J_{\Gamma_l,n}:D^2_{\xi_0,\eps}\times(D^2_{0,\eps}\setminus\{0\})\to(\CC\setminus\{0\})^2$,
$1\le n\le n_{\Gamma_l}$, таких что
$$
f\circ J_{\Gamma_l,n}(\xi,u)=\xi, \quad J_{\Gamma_l,n}^*(dz\wedge
dw)=\varkappa_{\Gamma_l,n}\ u^{(u_0-1)\alpha_{\Gamma_l}+(v_0-1)\beta_{\Gamma_l}-1}\ d\xi\wedge
du,
$$
$(\xi,u)\in D^2_{\xi_0,\eps}\times(D^2_{0,\eps}\setminus\{0\})$,
причем $\lim_{u\to0}|J_{\Gamma_l,n}(\xi,u)|=\infty$ равномерно по
$\xi\in D^2_{\xi_0,\eps}$, $1\le n\le n_{\Gamma_l}$, где
$n_{\Gamma_l}+1$ равно количеству целочисленных точек на стороне
$\Gamma_l$, $(\alpha_{\Gamma_l},\beta_{\Gamma_l})$ --- несократимый
вектор внешней нормали стороны $\Gamma_l$, $(u_0,v_0)\in\Gamma_l$
--- любая точка на стороне $\Gamma_l$, $\varkappa_{\Gamma_l,n}:=1$ при $(1,1)\not\in\Gamma_l$, $\varkappa_{\Gamma_l,(3\pm1)/2}:=\pm(a_{1,1}^2-4a_{2,0}a_{0,2})^{-1/2}$ при $(1,1)\in\Gamma_l$;

2) образы всех этих $n_\mu=\sum_{l}n_{\Gamma_l}$ вложений (отвечающих
одной и той же стороне, но разным значениям $n$, либо разным сторонам
многоугольника Ньютона) попарно не пересекаются, и
объединение этих образов содержит $f^{-1}(D^2_{\xi_0,\eps})\setminus
D^4_{0,R}$ (т.е.\ дополнение этого объединения в
$f^{-1}(D^2_{\xi_0,\eps})$ ограничено, а потому имеет компактное
замыкание в $\CC^2$);

3) каждое отображение
$D^2_{\xi_0,\eps}\times(D^2_{0,\eps}\setminus\{0\})\to(\CC\setminus\{0\})^2$, \\
$(\xi,u)\mapsto(u^{\alpha_{\Gamma_l}}\Pr_z(J_{\Gamma_l,n}(\xi,u)),u^{\beta_{\Gamma_l}}\Pr_w(J_{\Gamma_l,n}(\xi,u)))$,
продолжается до голоморфного отображения \\
$\hat J_{\Gamma_l,n}:D^2_{\xi_0,\eps}\times
D^2_{0,\eps}\to(\CC\setminus\{0\})^2$, такого что
$f^{(\alpha_{\Gamma_l},\beta_{\Gamma_l})}(\hat
J_{\Gamma_l,n}((\xi,0)))=0$, где $\Pr_z,\Pr_w:\CC^2\to\CC$ ---
проекции на первую и вторую компоненту соответственно, $f^{(\alpha_{\Gamma_l},\beta_{\Gamma_l})}(z,w)$ -- усеченный многочлен {\rm(см.\ определение~\ref {newton:def:0.1})}.
\end{theorem}

\begin{coment}
В действительности, теорема~\ref {newton:intrth2} и ее доказательство
(а также следствие~\ref {newton:cor2} без формул для
$\Delta_\xi|_{\SSigma_\xi\cap U_{\Gamma_l,n}}$ и
$g_\xi|_{\SSigma_\xi\cap U_{\Gamma_l,n}}$) останутся верными, если
вместо условия невырожденности многочлена $f$ относительно своего
многоугольника Ньютона (см.\ определение~\ref {def:nondegen}) и условий (i) и
$\dim\NP_{f-\xi_0}=2$ на многоугольник Ньютона $\NP_{f-\xi_0}$
наложить следующие (более слабые) условия.
Во-первых, потребовать выполнение условий определения~\ref
{def:nondegen} для любого ковектора $\eta\in\RR^{2*}$, {\it имеющего хотя бы одну
положительную координату}. Во-вторых, наложить условия (i) и
$\dim\NP_{f-\xi_0}=2$ не на многоугольник Ньютона $\NP_{f-\xi_0}$, а
на многоугольник
$\NP_{f-\xi_0}^*:=\conv\{(0,0)\}\cup(\cup_l\Gamma_l)$, где
$\cup_l\Gamma_l$ --- объединение всех сторон $\Gamma_l$
многоугольника Ньютона $\NP_{f-\xi_0}$, не лежащих на координатных
осях. В этом случае $\NP_{f-\xi_0}\subseteq\NP_{f-\xi_0}^*$.
\end{coment}

\begin{corollary} [О компактификации замкнутой окрестности слоя] \label {newton:cor2}
Пусть выполнены условия теоремы~$\ref {newton:intrth2}$. Имеется
комплексное 2-мерное связное многообразие $\widetilde M^4=\widetilde
M^4_{\xi_0,\eps}$ с комплексно аналитическим атласом из $n_\mu+1$
карт, полученное из
$M^4=M^4_{\xi_0,\eps}:=f^{-1}(D^2_{\xi_0,\eps})\subset\CC^2$
приклеиванием $n_\mu$ экземпляров множества $D^2_{\xi_0,\eps}\times
D^2_{0,\eps}\subset\CC^2$ при помощи вложений $J_{\Gamma_l,n}$
{\rm(см.\ теорему~\ref {newton:intrth2})}, такое что $\widetilde
M^4\setminus M^4\approx
D^2_{\xi_0,\eps}\times\{0\}\times\{1,\dots,n_\mu\}$ (``бесконечно
удаленные'' точки $p_{\xi,\Gamma_l,n}$) и замыкание
$\widetilde\SSigma_\xi\subset\widetilde M^4$ каждого слоя
$\SSigma_\xi=f^{-1}(\xi)$ в $\widetilde M^4$, является компактным
слоем $\widetilde\SSigma_\xi=\widetilde f^{-1}(\xi)$ некоторой
голоморфной функции $\widetilde f:\widetilde M^4\to\CC$, $\xi\in
D^2_{\xi_0,\eps}$, такой что $\widetilde f|_{M^4}=f$ и множества
критических точек функций $f$ и $\widetilde f$ совпадают. В
$U_{\Gamma_l,n}:=J_{\Gamma_l,n}(D^2_{\xi_0,\eps}\times(D^2_{0,\eps}\setminus\{0\}))\subset
M^4$ векторное поле $\sgrad_\CC f$, 1-форма $\Delta_\xi$ и риманова
метрика пополнения $g_\xi$ {\rm(см.\ определения~\ref
{newton:def:0.4} и~\ref {newton:def:0.6}(А))} имеют следующий вид в
координатах $(\xi,u)\in
D^2_{\xi_0,\eps}\times(D^2_{0,\eps}\setminus\{0\})$ из теоремы~$\ref
{newton:intrth2}$:
$$
\sgrad_\CC
f|_{U_{\Gamma_l,n}}=\frac1{\varkappa_{\Gamma_l,n}} u^{(1-u_0)\alpha_{\Gamma_l}+(1-v_0)\beta_{\Gamma_l}+1}\frac{\d}{\d
u},
$$
$$
\Delta_\xi|_{\SSigma_\xi\cap
U_{\Gamma_l,n}}=\varkappa_{\Gamma_l,n}\, u^{(u_0-1)\alpha_{\Gamma_l}+(v_0-1)\beta_{\Gamma_l}-1}du,
\qquad \xi\in D^2_{\xi_0,\eps},
$$
$$
g_\xi|_{\SSigma_\xi\cap U_{\Gamma_l,n}}=|\varkappa_{\Gamma_l,n}|^2(u\overline
u)^{(u_0-1)\alpha_{\Gamma_l}+(v_0-1)\beta_{\Gamma_l}-1}du\,d\overline
u, \qquad \xi\in D^2_{\xi_0,\eps}.
$$
При этом $(u_0-1)\alpha_{\Gamma_l}+(v_0-1)\beta_{\Gamma_l}-1\ge0$
тогда и только тогда, когда многоугольник Ньютона $\NP_{f-\xi_0}$
содержит хотя бы одну внутреннюю точку с целыми координатами (т.е.\
когда $\widetilde\SSigma_{\xi_0}\not\approx S^2$).
\end{corollary}

\begin{coment}
При $|\xi-\xi_0|<\eps$ векторное поле $\sgrad_\CC
f|_{\widetilde\SSigma_\xi}$ на компактной связной поверхности
$\widetilde\SSigma_\xi=\widetilde f^{-1}(\xi)$ имеет ровно $n_\mu$
особых точек, индексы которых равны
$(1-u_0)\alpha_{\Gamma_l}+(1-v_0)\beta_{\Gamma_l}+1$, см.\
следствие~\ref {newton:cor2}. Поэтому сумма индексов равна
$n_\mu-2S(\NP_{f-\xi})$, где $S(\NP_{f-\xi})$ --- площадь
многоугольника $\conv\{(1,1)\}\cup(\cup_l\Gamma_l)$. С другой
стороны, по следствию~\ref {newton:cor:1.3} род поверхности
$\widetilde\SSigma_\xi$ (т.е.\ количество ручек) равен
$n_g=B^+(\NP_{f-\xi})$. Так как сумма индексов особых точек
векторного поля равна $2-2n_g$, получаем равенство
$n_\mu-2S(\NP_{f-\xi})=2-2B^+(\NP_{f-\xi})$, равносильное известной
теореме Пика.
\end{coment}

\begin{proof}[Доказательство теоремы \ref {newton:intrth2}.]
Сначала отметим, что в силу условия (i) выполнено
$\alpha_{\Gamma_l}\ge0$ и $\beta_{\Gamma_l}\ge0$, причем по крайней
мере одно неравенство строгое. Пусть для определенности $(u_0,v_0)$
--- начальная вершина стороны $\Gamma_l\subset\d\NP_{f-\xi_0}$ по
отношению к положительной ориентации (против часовой стрелки)
замкнутой ломаной $\d\NP_{f-\xi_0}\subset\CC$. Тогда целочисленные
точки стороны $\Gamma_l$ имеют координаты
$(u_n,v_n):=(u_0-n\beta_{\Gamma_l},v_0+n\alpha_{\Gamma_l})$,
$n=0,1,\dots,n_{\Gamma_l}$. Далее, рассмотрим вектор
$\eta:=(\alpha_{\Gamma_l},\beta_{\Gamma_l})$ внешней нормали стороны
$\Gamma_l$ и отвечающий ему усеченный многочлен $f^{\eta}(z,w)
 =\sum\limits_{n=0}^{n_{\Gamma_l}}a_{u_n,v_n}z^{u_n}w^{v_n}
 =z^{u_0}w^{v_0} P_{\Gamma_l}\left(\frac{w^{\alpha_{\Gamma_l}}}{z^{\beta_{\Gamma_l}}}\right)$, где
$$
P_{\Gamma_l}(y):=\sum\limits_{n=0}^{n_{\Gamma_l}}a_{u_n,v_n}y^n=\sum\limits_{n=0}^{n_{\Gamma_l}}a_{u_0-n\beta_{\Gamma_l},v_0+n\alpha_{\Gamma_l}}y^n.
$$
Пусть $\hat{y}_1,\dots,\hat{y}_{n_{\Gamma_l}}$ --- корни уравнения
$P_{\Gamma_l}(y)=0$. Заметим, что все корни уравнения
$P_{\Gamma_l}(y)=0$ различны (так как $f(z,w)$ --- невырожденный
многочлен) и не равны нулю (так как $a_{u_0,v_0}\ne0$ и
$a_{u_{n_{\Gamma_l}},v_{n_{\Gamma_l}}}=a_{u_0-n_{\Gamma_l}\beta_{\Gamma_l},v_0+n_{\Gamma_l}\alpha_{\Gamma_l}}\ne0$
в силу того, что $(u_0,v_0)$ и
$(u_{n_{\Gamma_l}},v_{n_{\Gamma_l}})=(u_0-n_{\Gamma_l}\beta_{\Gamma_l},v_0+n_{\Gamma_l}\alpha_{\Gamma_l})$
--- вершины многоугольника Ньютона). Положим $\hat{y}:=\hat{y}_n$,
где $n=1,\dots,n_{\Gamma_l}$.

Шаг 1. Пусть $\alpha_{\Gamma_l}\ne0$. Рассмотрим отображение
$$
I_{\eps_2,\eps_3,\alpha_{\Gamma_l},\beta_{\Gamma_l},\hat{y}}:(D^2_{0,\eps_2}\setminus\{0\})\times
D^2_{0,\eps_3}\to\CC^2, \quad
(u,g)\mapsto(u^{-\alpha_{\Gamma_l}},u^{-\beta_{\Gamma_l}}(\hat{y}+g)^{1/\alpha_{\Gamma_l}}),
$$
где $\eps_2,\eps_3>0$ --- некоторые числа.
Покажем, что существует единственная непрерывная функция
$g_{\Gamma_l,n}(\xi,u)$, такая что выполнены соотношения
$$
g_{\Gamma_l,n}(\xi,0)=0, \quad f\circ
I_{\eps_2,\eps_3,\alpha_{\Gamma_l},\beta_{\Gamma_l},\hat{y}}(u,g_{\Gamma_l,n}(\xi,u))=\xi.
$$
В самом деле, значение любого монома $a_{k,m}z^kw^m$ многочлена
$f(z,w)$ на паре
$(z,w)=I_{\eps_2,\eps_3,\alpha_{\Gamma_l},\beta_{\Gamma_l},\hat{y}}(u,g)$
равно
$a_{k,m}z^kw^m=a_{k,m}u^{-k\alpha_{\Gamma_l}-m\beta_{\Gamma_l}}(\hat{y}+g)^{m/\alpha_{\Gamma_l}}
=a_{k,m}u^{-\alpha_{\Gamma_l}u_0-\beta_{\Gamma_l}v_0}u^{\alpha_{\Gamma_l}(u_0-k)+\beta_{\Gamma_l}(v_0-m)}
(\hat{y}+g)^{m/\alpha_{\Gamma_l}}$. Заметим, что
$\alpha_{\Gamma_l}(u_0-k)+\beta_{\Gamma_l}(v_0-m)\ge0$, поскольку
$(k,m)\in\NP_{f-\xi_0}$, и
$\alpha_{\Gamma_l}(u_0-k)+\beta_{\Gamma_l}(v_0-m)=0$ тогда и только
тогда, когда $(k,m)\in\Gamma_l$, поскольку
$(\alpha_{\Gamma_l},\beta_{\Gamma_l})$ --- внешняя нормаль к
$\Gamma_l$. Отсюда следует, что $f\circ
I_{\eps_2,\eps_3,\alpha_{\Gamma_l},\beta_{\Gamma_l},\hat{y}}(u,g)-\xi
 =u^{-\alpha_{\Gamma_l}u_0-\beta_{\Gamma_l}v_0}(P_{\Gamma_l}(\hat{y}+g)+\dots)$,
где невыписанные мономы имеют степень по переменной $u$ больше либо
равную единице. Рассмотрим в области $\CC\times D^2_{0,\eps_2}\times
D^2_{0,\eps_3}$ голоморфную функцию
$$
F(\xi,u,g):=(f\circ
I_{\eps_2,\eps_3,\alpha_{\Gamma_l},\beta_{\Gamma_l},\hat{y}}(u,g)-\xi)u^{\alpha_{\Gamma_l}u_0+\beta_{\Gamma_l}v_0}=P_{\Gamma_l}(\hat{y}+g)+\dots
$$
и рассмотрим уравнение $F(\xi,u,g)=0$. Заметим, что для любого $\xi$
при $u=0$ подстановка $g=0$ дает решение, поскольку $\hat{y}$ ---
корень уравнения $P_{\Gamma_l}(y)=0$, т.е. $F(\xi,0,0)=0$. Далее,
$(\partial F/\partial
g)|_{(\xi,0,0)}=(P'_{\Gamma_l}(\hat{y}+g)+\dots)|_{(\xi,0,0)}=P'_{\Gamma_l}(\hat{y})\ne0$,
так как (в силу невырожденности многочлена $f-\xi_0$ относительно
своего многогоугольника Ньютона $\NP_{f-\xi_0}$) многочлен
$P_{\Gamma_l}(y)$ не имеет кратных корней. Поэтому, по теореме о
неявной функции, существуют $\eps_1,\eps_2,\eps_3>0$, такие что
существует единственная функция
$g_{\Gamma_l,n}:D^2_{\xi_0,\eps_1}\times D^2_{0,\eps_2}\to
D^2_{0,\eps_3}\subset\CC$ со свойством
 $$
g_{\Gamma_l,n}(\xi,0)=0, \quad F(\xi,u,g_{\Gamma_l,n}(\xi,u))=0,
 $$
а, стало быть,
$I_{\eps_2,\eps_3,\alpha_{\Gamma_l},\beta_{\Gamma_l},\hat{y}}(u,g_{\Gamma_l,n}(\xi,u))\in\SSigma_\xi$.
Более того, согласно теореме о неявной функции, функция
$g_{\Gamma_l,n}=g_{\Gamma_l,n}(\xi,u)$ является голоморфной.

Шаг 2. Далее, положим
 $$
  J_{\eps_1,\eps_2,\eps_3,\alpha_{\Gamma_l},\beta_{\Gamma_l},\hat{y}}(\xi,u):=
I_{\eps_2,\eps_3,\alpha_{\Gamma_l},\beta_{\Gamma_l},\hat{y}}(u,g_{\Gamma_l,n}(\xi,u)).
 $$
Отсюда и из определения функции $F$ и свойства
$F(\xi,u,g_{\Gamma_l,n}(\xi,u))=0$ получаем \\
$f\circ J_{\eps_1,\eps_2,\eps_3,\alpha_{\Gamma_l},\beta_{\Gamma_l},\hat{y}}(\xi,u)=\xi$.
Вычислим
$$
J_{\eps_1,\eps_2,\eps_3,\alpha_{\Gamma_l},\beta_{\Gamma_l},\hat{y}}^{*}(dz\wedge
dw)
 =d(u^{-\alpha_{\Gamma_l}})\wedge d(u^{-\beta_{\Gamma_l}}(\hat{y}+g_{\Gamma_l,n}(\xi,u))^{1/\alpha_{\Gamma_l}})
$$
$$
 = (-\alpha_{\Gamma_l}u^{-\alpha_{\Gamma_l}-1}du)\wedge(
\frac{(\hat{y}+g_{\Gamma_l,n}(\xi,u))^{(1-\alpha_{\Gamma_l})/\alpha_{\Gamma_l}}}{\alpha_{\Gamma_l}}
u^{-\beta_{\Gamma_l}}\frac{\partial
g_{\Gamma_l,n}(\xi,u)}{\partial\xi}d\xi)
 $$
 $$
=
(\hat{y}+g_{\Gamma_l,n}(\xi,u))^{(1-\alpha_{\Gamma_l})/\alpha_{\Gamma_l}}u^{-\alpha_{\Gamma_l}-\beta_{\Gamma_l}-1}
(-\frac{\partial F(\xi,u,g_{\Gamma_l,n})}{\partial\xi}/\frac{\partial
F(\xi,u,g_{\Gamma_l,n})}{\partial g})d\xi\wedge du
 $$
 $$
=u^{\alpha_{\Gamma_l}(u_0-1)+\beta_{\Gamma_l}(v_0-1)-1}\frac{(\hat{y}+g_{\Gamma_l,n}(\xi,u))^{(1-\alpha_{\Gamma_l})/\alpha_{\Gamma_l}}}
{P'_{\Gamma_l}(\hat{y}+g_{\Gamma_l,n}(\xi,u))+\dots}d\xi\wedge du,
$$
где невыписанные мономы, как и прежде, имеют степень по $u$ больше
либо равную единице. Поэтому
 $$
 J_{\eps_1,\eps_2,\eps_3,\alpha_{\Gamma_l},\beta_{\Gamma_l},\hat{y}}^{*}(dz\wedge dw)=
u^{\alpha_{\Gamma_l}(u_0-1)+\beta_{\Gamma_l}(v_0-1)-1}h(\xi,u)d\xi\wedge
du,
 $$
где $h(\xi,u)$ --- голоморфная отделенная от нуля функция двух
комплексных переменных в некоторой окрестности точки $(\xi_0,0)$.
Отсюда существуют $\eps>0$ и замена координат $\psi=\psi_{\Gamma_l,n}:U(\xi_0,0)\to
D^2_{\xi_0,\eps}\times D^2_{0,\eps}$, $\psi:(\xi,u)\mapsto(\tilde{\xi},\tilde{u})$, такие что $\tilde{\xi}\equiv\xi$ и
 $$
(J_{\eps_1,\eps_2,\eps_3,\alpha_{\Gamma_l},\beta_{\Gamma_l},\hat{y}}\circ\psi^{-1})^{*}(dz\wedge
dw)=\varkappa_{\Gamma_l,n}\
\tilde{u}^{\alpha_{\Gamma_l}(u_0-1)+\beta_{\Gamma_l}(v_0-1)-1}d\xi\wedge
d\tilde{u},
 $$
где $U(\xi_0,0)\subset\CC^2$ --- некоторая окрестность точки $(\xi_0,0)$ в $\CC^2$, $\varkappa_{\Gamma_l,n}\in\CC\setminus\{0\}$ как в формулировке теоремы.
 Положим $J_{\Gamma_l,n}:=
 J_{\eps_1,\eps_2,\eps_3,\alpha_{\Gamma_l},\beta_{\Gamma_l},\hat{y}}\circ\psi^{-1}$.

Пункт 3) следует из того, что $\hat J_{\Gamma_l,n}(\xi,u)=(1,\hat
y+g_n(\xi,u))$, $g_n(\xi,0)=0$, $\hat y\ne0$ и
$f^{(\alpha_{\Gamma_l},\beta_{\Gamma_l})}(1,\hat y)=P_{\Gamma_l}(\hat
y)=0$. Из пункта 3) следует, что
$\lim_{u\to0}|J_{\Gamma_l,n}(\xi,u)|=\infty$ равномерно по $\xi\in
D^2_{\xi_0,\eps}$.

Шаг 3. Докажем, что
$J_{\eps_1,\eps_2,\eps_3,\alpha_{\Gamma_l},\beta_{\Gamma_l},\hat{y}}$
--- вложение.
Допустим противное, тогда существуют $(\xi_1,u_1)$ и $(\xi_2,u_2)$,
такие что $(\xi_1,u_1)\ne(\xi_2,u_2)$ и
$J_{\eps_1,\eps_2,\eps_3,\alpha_{\Gamma_l},\beta_{\Gamma_l},\hat{y}}(\xi_1,u_1)=
J_{\eps_1,\eps_2,\eps_3,\alpha_{\Gamma_l},\beta_{\Gamma_l},\hat{y}}(\xi_2,u_2)$.
Отсюда $\xi_1=\xi_2$, так как
$\xi_1=f(J_{\eps_1,\eps_2,\eps_3,\alpha_{\Gamma_l},\beta_{\Gamma_l},\hat{y}}(\xi_1,u_1))=
f(J_{\eps_1,\eps_2,\eps_3,\alpha_{\Gamma_l},\beta_{\Gamma_l},\hat{y}}(\xi_2,u_2))=\xi_2$,
обозначим $\xi:=\xi_1=\xi_2$. С другой стороны, так как
$(u_1^{-\alpha_{\Gamma_l}},u_1^{-\beta_{\Gamma_l}}(\hat{y}+g_{\Gamma_l,n}(\xi,u_1))^{1/\alpha_{\Gamma_l}})=
(u_2^{-\alpha_{\Gamma_l}},u_2^{-\beta_{\Gamma_l}}(\hat{y}+g_{\Gamma_l,n}(\xi,u_2))^{1/\alpha_{\Gamma_l}})$,
то $g_{\Gamma_l,n}(\xi,u_1)=g_{\Gamma_l,n}(\xi,u_2)$ и (в силу
взаимной простоты $\alpha_{\Gamma_l}$ и $\beta_{\Gamma_l}$)
$u_1=u_2$, противоречие.

Аналогично определяется вложение
$J_{\eps_1,\eps_2,\eps_3,\alpha_{\Gamma_l},\beta_{\Gamma_l},\hat{y}}$,
если $\alpha_{\Gamma_l}=0$ и $\beta_{\Gamma_l}\ne0$ (как следствие,
$\beta_{\Gamma_l}=1$). Тем самым, пункты 1) и 3) доказаны. Докажем
пункт 2).

Шаг 4. Докажем, что образы построенных вложений попарно не
пересекаются. Пусть $\Gamma_{l_1}\ne\Gamma_{l_2}$. Возможны два
случая.

{\it Случай 1.} Пусть $\alpha_{\Gamma_{l_i}}\ne0$, $i=1,2$. Допустим,
что существуют $n_1,n_2\in\NN$, $1\le n_1\le n_{\Gamma_{l_1}}$, $1\le n_2\le
n_{\Gamma_{l_2}}$, такие что для любых сколь угодно малых
$\eps_1,\eps_2>0$ существуют $(\xi_i,u_i)$, $|\xi_i-\xi_0|<\eps_1$,
$|u_i|<\eps_2$, $i=1,2$, такие что выполнено соотношение: $
J_{\eps_1,\eps_2,\eps_3,\alpha_{\Gamma_{l_1}},\beta_{\Gamma_{l_1}},\hat{y}_{n_1}}(\xi_1,u_1)=
J_{\eps_1,\eps_2,\eps_3,\alpha_{\Gamma_{l_2}},\beta_{\Gamma_{l_2}},\hat{y}_{n_2}}(\xi_2,u_2).
$ Тогда
 $$
u_1^{-\alpha_{\Gamma_{l_1}}}=u_2^{-\alpha_{\Gamma_{l_2}}}, \quad
u_1^{-\beta_{\Gamma_{l_1}}}(\hat{y}_{n_1}+g_{\Gamma_{l_1},n_1})^{1/\alpha_{\Gamma_{l_1}}}=
u_2^{-\beta_{\Gamma_{l_2}}}(\hat{y}_{n_2}+g_{\Gamma_{l_2},n_2})^{1/\alpha_{\Gamma_{l_2}}}.
 $$
Отсюда $u_2^{\alpha_{\Gamma_{l_1}}\beta_{\Gamma_{l_2}}-
\alpha_{\Gamma_{l_2}}\beta_{\Gamma_{l_1}}}=\frac
{(\hat{y}_{n_2}+g_{\Gamma_{l_2},n_2})^{\alpha_{\Gamma_{l_1}}/\alpha_{\Gamma_{l_2}}}}
{\hat{y}_{n_1}+g_{\Gamma_{l_1},n_1}}$, а, значит,
$u_2^{\alpha_{\Gamma_{l_1}}\beta_{\Gamma_{l_2}}-
\alpha_{\Gamma_{l_2}}\beta_{\Gamma_{l_1}}}$ ограничено и отделено от
нуля (при сколь угодно малых $u_2$). Следовательно,
$\alpha_{\Gamma_{l_1}}\beta_{\Gamma_{l_2}}-
\alpha_{\Gamma_{l_2}}\beta_{\Gamma_{l_1}}=0$, откуда векторы
$(\alpha_{\Gamma_{l_1}},\beta_{\Gamma_{l_1}})$ и
$(\alpha_{\Gamma_{l_2}},\beta_{\Gamma_{l_2}})$ пропорциональны, а
потому совпадают, противоречие.

{\it Случай 2.} Пусть $\alpha_{\Gamma_{l_1}}\ne0$ и
$\alpha_{\Gamma_{l_2}}=0$ (заметим, что одновременно равенства нулю
$\alpha_{\Gamma_{l_i}}=0$, $i=1,2$, невозможны, в силу условия
$\Gamma_{l_1}\ne\Gamma_{l_2}$). Допустим, что существуют
$n_1,n_2\in\NN$, $1\le n_1\le n_{\Gamma_{l_1}}$, $1\le n_2\le n_{\Gamma_{l_2}}$, такие что
для любых сколь угодно малых $\eps_1,\eps_2>0$ существуют
$(\xi_i,u_i)$, $|\xi_i-\xi_0|<\eps_1$, $|u_i|<\eps_2$, $i=1,2$, такие
что выполнено соотношение: $
J_{\eps_1,\eps_2,\eps_3,\alpha_{\Gamma_{l_1}},\beta_{\Gamma_{l_1}},\hat{y}_{n_1}}(\xi_1,u_1)=
J_{\eps_1,\eps_2,\eps_3,0,1,\hat{y}_{n_2}}(\xi_2,u_2). $ Тогда
$u_1^{-\alpha_{\Gamma_{l_1}}}=(\hat{y}_{n_2}+g_{\Gamma_{l_2},n_2})^{1/
\beta_{\Gamma_{l_2}}}$. Отсюда $u_1^{-\alpha_{\Gamma_{l_1}}}$
ограничено и отделено от нуля (при сколь угодно малых $u_1$), что
невозможно в силу условия $\alpha_{\Gamma_{l_1}}>0$, противоречие.

Далее, рассмотрим сторону $\Gamma_{l}$. Пусть
$\alpha_{\Gamma_{l}}\ne0$. Допустим, что существуют
$n_1,n_2\in\NN$,  $1\le n_1<n_2\le n_{\Gamma_{l_1}}$, такие что для
любых сколь угодно малых $\eps_1,\eps_2>0$ существуют $(\xi_i,u_i)$,
$|\xi_i-\xi_0|<\eps_1$, $|u_i|<\eps_2$, $i=1,2$, такие что выполнено
соотношение: $
J_{\eps_1,\eps_2,\eps_3,\alpha_{\Gamma_{l}},\beta_{\Gamma_{l}},\hat{y}_{n_1}}(\xi_1,u_1)=
J_{\eps_1,\eps_2,\eps_3,\alpha_{\Gamma_{l}},\beta_{\Gamma_{l}},\hat{y}_{n_2}}(\xi_2,u_2).
$ Тогда
$\hat{y}_{n_1}-\hat{y}_{n_2}=g_{\Gamma_l,n_2}(\xi_2,u_2)-g_{\Gamma_l,n_1}(\xi_1,u_1)$,
что невозможно, в силу того, что $\hat{y}_{n_1}\ne\hat{y}_{n_2}$ и
$g_{\Gamma_l,n}(\xi,u)$ --- ограниченная голоморфная функция, причем
$g_{\Gamma_{l,n}}(\xi,0)=0$, противоречие. Аналогично рассматривается
случай, когда $\beta_{\Gamma_l}\ne0$.

Шаг 5. Докажем, что при $0<\eps_1<|a_{0,0}-\xi_0|/2$ множество
 $$
f^{-1}(D^2_{\xi_0,\eps_1})\setminus\bigcup\limits_{\Gamma_l,n}J_{\eps_1,\eps_2,\eps_3,\alpha_{\Gamma_{l}},\beta_{\Gamma_{l}},\hat{y}_n}
\left(D^2_{\xi_0,\eps_1}\times (D^2_{0,\eps_2}\setminus\{0\})\right)
 $$
ограничено в $\CC^2$. Обозначим
$X:=\bigcup\limits_{\Gamma_l,n}J_{\eps_1,\eps_2,\eps_3,\alpha_{\Gamma_{l}},\beta_{\Gamma_{l}},\hat{y}_n}
(D^2_{\xi_0,\eps_1}\times (D^2_{0,\eps_2}\setminus\{0\}))$.
Осталось показать, что $f^{-1}(D^2_{\xi_0,\eps_1})\setminus X$
ограничено. Допустим противное, тогда существует последовательность
$(z_j, w_j)\in f^{-1}(D^2_{\xi_0,\eps_1})\setminus X$, $j\in\NN$,
такая что либо $z_j\to\infty$, либо $w_j\to\infty$. Тогда
$|z_j|=e^{\alpha_j}, |w_j|=e^{\beta_j}$, где $\alpha_j,
\beta_j\in\RR$ и $\max\{\alpha_j,\beta_j\}\to+\infty$.
Возможны два случая.

{\it Случай 1.} Допустим, что для любой стороны $\Gamma_{l}$
многоугольника Ньютона последовательность $(\alpha_j, \beta_j)$
отделена с точностью до пропорциональности от вектора нормали
$(\alpha_{\Gamma_{l}}, \beta_{\Gamma_{l}})$.
Тогда $|f(z_j,w_j)-\xi_j| = e^{u_0\alpha_j + v_0\beta_j}(|\tilde
a_{u_0,v_0}|+o(1))$ при $j\to\infty$, где $\xi_j:=f(z_j,w_j)$,
$\tilde a_{u,v}:=a_{u,v}$ при $(u,v)\ne(0,0)$, $\tilde
a_{0,0}:=a_{0,0}-\xi_j$, $u_0,v_0\in\ZZ$ --- координаты вершины
многоугольника Ньютона, на которой значение выражения $\alpha_j u_0 +
\beta_j v_0$ наибольшее (такая вершина в рассматриваемом случае, без
ограничения общности, не зависит от $j$ при достаточно большом $j$).
Левая часть полученного равенства равна нулю, а правая отлична от
нуля, противоречие.

{\it Случай 2.} Таким образом, последовательность $(\alpha_j,
\beta_j)$ с точностью до пропорциональности стремится к вектору
внешней нормали $(\alpha_{\Gamma_{l}}, \beta_{\Gamma_{l}})$ некоторой
стороны $\Gamma_{l}$ многоугольника Ньютона. Возможны два подслучая.

{\it Подслучай 2а.} Допустим, что $(\alpha_{\Gamma_{l}},
\beta_{\Gamma_{l}})=(-1,0)$.
Тогда $|z_j|=e^{-t_j}\to0$, $|w_j|=e^{o(t_j)}\to+\infty$, где
$t_j\to+\infty$ при $j\to\infty$, и для любого монома
$|z_j^pw_j^q|=e^{t_j(-p+o(1))}$, $p,q\in\ZZ$. Отсюда следует, что
$0=f(z_j,w_j)-\xi_j=o(1)+f(0,w_j)-\xi_j$, где $\xi_j:=f(z_j,w_j)$.
Поэтому, ввиду того, что
$|w_j|\to\infty$, имеем $f(0,w)\equiv\const$, то есть
$f(z,w)-\xi=zL(z,w)$ для некоторых константы $\xi\in\CC$ и многочлена
$L=L(z,w)$. Получили противоречие ввиду условий (i) и $\dim\NP_{f-\xi_0}=2$.

{\it Подслучай 2б.} Таким образом, $\alpha_{\Gamma_l}\ge0$ и
$\beta_{\Gamma_l}\ge0$. Напомним, что целочисленные точки стороны
$\Gamma_l$ имеют координаты
$(u_0-n\beta_{\Gamma_l},v_0+n\alpha_{\Gamma_l})$,
$n=0,1,\dots,n_{\Gamma_l}$, где $(u_0,v_0)$ --- начальная вершина
стороны $\Gamma_l\subset\d\NP_{f-\xi_0}$ по отношению к
положительной ориентации (против часовой стрелки) замкнутой ломаной
$\d\NP_{f-\xi_0}\subset\CC$. Пусть для определенности
$\alpha_{\Gamma_l}\ne0$ (тогда $\alpha_{\Gamma_l}>0$,
$\beta_{\Gamma_l}\ge0$ и $|z_j|\to\infty$), рассмотрим
($\alpha_{\Gamma_l}$-значное) отображение
 $$
h_{\Gamma_l}:(x,y)\mapsto(z,w)=(x^{\alpha_{\Gamma_l}},x^{\beta_{\Gamma_l}}y^{1/\alpha_{\Gamma_l}})
 $$
области $(\CC\setminus\{0\})\times\CC$ в себя. Тогда
$z_j=x_j^{\alpha_{\Gamma_l}}$,
$w_j=x_j^{\beta_{\Gamma_l}}y_j^{1/\alpha_{\Gamma_l}}$ для некоторых
$x_j,y_j\in\CC$. Имеем $x_j\to\infty$, $|y_j|=|x_j|^{o(1)}$, откуда
 $$
|f(z_j,w_j)-f_{\Gamma_l}(z_j,w_j)-\xi_j| =
O(|x_j|^{u_0\alpha_{\Gamma_l} + v_0\beta_{\Gamma_l}-1/2}),
 $$
$$
f_{\Gamma_l}(z_j,w_j)
 = x_j^{u_0\alpha_{\Gamma_l}+v_0\beta_{\Gamma_l}}
   \left(\sum_{n=0}^{n_{\Gamma_l}}a_{u_n,v_n} y_j^{v_n/\alpha_{\Gamma_l}}\right)
 = x_j^{u_0\alpha_{\Gamma_l}+v_0\beta_{\Gamma_l}} y_j^{v_0/\alpha_{\Gamma_l}}
   \left(\sum_{n=0}^{n_{\Gamma_l}}a_{u_n,v_n} y_j^n\right),
$$
где $u_n:=u_0-n\beta_{\Gamma_l}$, $v_n:=v_0+n\alpha_{\Gamma_l}$.
Поэтому
$$
0=f(z_j,w_j)-\xi_j=x_j^{u_0\alpha_{\Gamma_l}+v_0\beta_{\Gamma_l}}
   \left(O(|x_j|^{-1/2})+y_j^{v_0/\alpha_{\Gamma_l}}\sum_{n=0}^{n_{\Gamma_l}}a_{u_n,v_n}  y_j^n\right)
$$
при $j\to\infty$. Отсюда получаем, что последовательность
$y_j\in\CC$ может стремиться только к бесконечности (при $v_0<0$,
$v_{n_{\Gamma_l}}<0$), к нулю (при $v_0>0$, $v_{n_{\Gamma_l}}>0$) и к
корням многочлена
$P_{\Gamma_l}(y):=\sum_{n=0}^{n_{\Gamma_l}}a_{u_n,v_n}y^n$.
Первое невозможно, так как $f(z,w)$ является обычным многочленом (а
не многочленом Лорана), а потому $v_0\ge0$ и $v_{n_{\Gamma_l}}\ge0$.
Второе тоже невозможно, так как в противном случае выполнялось бы
$|x_j|\to\infty$ и $|y_j|=|x_j|^{-\eps_j}\to0$ для некоторого
$\eps_j\to0$, откуда $\eps_j>0$ (начиная с некоторого $j$), и
$|x_j|^{-1/2}=|y_j|^{1/(2\eps_j)}=o(|y_j|^{v_0/\alpha_{\Gamma_l}})$,
что приводит к противоречию.
Поэтому последовательность $y_j\in\CC$ ограничена
и может иметь своими предельными точками только корни многочлена
$P_{\Gamma_l}(y)$.

Без ограничения общности будем считать, что
$\lim_{j\to\infty}y_j=\hat y_n=:\hat y\in\CC$ и
$P_{\Gamma_l}(y)=(y-\hat y)Q(y)$, где $Q$ --- многочлен степени
$n_{\Gamma_l}-1$, тогда $\hat y\ne0$ в силу условия невырожденности.

Определим переменные $(u,g)=h_{\Gamma_l,\hat y}(x,y):=(x^{-1},y-\hat
y)$, тогда в этих переменных в окрестности начала координат отношение
$(f\circ
h_{\Gamma_l}-\xi)/x^{u_0\alpha_{\Gamma_l}+v_0\beta_{\Gamma_l}}$
является голоморфной функцией $F(\xi,u,g)$, совпадающей с введенной
на шаге 1. Так как $(u_j,g_j)\to(0,0)$ (в силу $x_j\to\infty$ и
$y_j\to\hat y$ по доказанному выше), $|\xi_j-\xi_0|<\eps_1$ и
$F(\xi_j,u_j,g_j)=0$, то в силу шага 1 имеем
$g_j=g_{\Gamma_l,n}(\xi_j,u_j)$, начиная с некоторого $j$, где
$\xi_j:=f(z_j,w_j)$ и $g_{\Gamma_l,n}(\xi,u)$ --- функция из шага 1.
Отсюда и из шага 2 имеем
$$
(z_j,w_j) =h_{\Gamma_l}\circ h_{\Gamma_l,\hat y}^{-1}(u_j,g_j)
=(u_j^{-\alpha_{\Gamma_l}},u_j^{-\beta_{\Gamma_l}}(\hat{y}+g_j)^{1/\alpha_{\Gamma_l}})
$$
$$
=I_{\eps_2,\eps_3,\alpha_{\Gamma_l},\beta_{\Gamma_l},\hat{y}}(u_j,g_j)
=J_{\eps_1,\eps_2,\eps_3,\alpha_{\Gamma_l},\beta_{\Gamma_l},\hat{y}}(\xi_j,u_j).
$$
Отсюда $(z_j,w_j)\in X$, начиная с некоторого $j$, противоречие.
\end{proof}

\begin{corollary} [Условие компактности пополненных слоев, типы особенностей поля $\sgrad f|_{\SSigma_\xi}$] \label {newton:cor:Newton}
Пусть $f(z,w)-\xi_0$ --- невырожденный многочлен относительно своего
многоугольника Ньютона $\NP_{f-\xi_0}$, удовлетворяющий условию {\rm
(i)}, причем $\dim\NP_{f-\xi_0}=2$, где
$\dim\NP_{f-\xi_0}$ --- размерность минимального линейного
пространства в $\RR^2$, содержащего $\NP_{f-\xi_0}$. Тогда существует
$\eps>0$ такое, что для любого $\xi\in\CC$, $|\xi-\xi_0|\le\eps$,
выполнено:

(А) многочлен $f(z,w)-\xi$ невырожден относительно своего
многоугольника Ньютона $\NP_{f-\xi}$, $\xi$ --- неособое значение,
$\dim\NP_{f-\xi}=2$;

(Б) при $n_g:=B^+(\NP_{f-\xi})\ge1$ поток векторного поля $\sgrad_\CC
f|_{\SSigma_\xi}$ является неполным, и пополнение
$\overline{\SSigma}_\xi$ слоя $\SSigma_\xi=f^{-1}(\xi)$ относительно
метрики пополнения $g_\xi$ является компактной связной поверхностью с
плоской метрикой и коническими особенностями; эта поверхность
гомеоморфна слою $\widetilde\SSigma_\xi=\widetilde f^{-1}(\xi)$
{\rm(см.\ следствие~\ref {newton:cor2})}, гомеоморфна сфере с $n_g$
ручками, причем $|\overline{\SSigma}_\xi\setminus\SSigma_\xi|=n_\mu$,
и в точках множества $\overline{\SSigma}_\xi\setminus\SSigma_\xi$,
называемых {\em бесконечно удаленными}, метрика имеет конические
особенности; при $n_g=0$ поток векторного поля $\sgrad_\CC
f|_{\SSigma_\xi}$ является полным, и пополнение
$\overline{\SSigma}_\xi$ любого слоя $\SSigma_\xi$ относительно
метрики пополнения $g_\xi$ совпадает с самим слоем $\SSigma_\xi$ и
изометрично евклидовой плоскости или плоскому цилиндру;

(В) различным сторонам многоугольника Ньютона, не лежащим на
координатных осях, соответствуют различные бесконечно удаленные точки
на $\overline{\SSigma}_\xi$; количество различных бесконечно
удаленных точек, отвечающих одной и той же стороне $\Gamma_l$
многоугольника Ньютона, равно $n_{\Gamma_l}$; при $n_g\ge1$
(соответственно $n_g=0$) в каждой из этих точек векторное поле
$\sgrad f|_{\SSigma_\xi}$ имеет особенность полюс порядка
$(u_0-1)\alpha_{\Gamma_l}+(v_0-1)\beta_{\Gamma_l}-1\ge0$
(соответственно ноль порядка
$(1-u_0)\alpha_{\Gamma_l}+(1-v_0)\beta_{\Gamma_l}+1>0$); при
$n_g\ge1$ плоская метрика на пополненном слое
$\overline{\SSigma}_\xi$ имеет в каждой бесконечно удаленной точке
коническую особенность с полным углом
$2\pi((u_0-1)\alpha_{\Gamma_l}+(v_0-1)\beta_{\Gamma_l})$.
\end{corollary}

\begin{proof}
Пункт (А) следует из теоремы~\ref{newton:lem:dis} и того факта, что в
силу условия (i) многоугольник Ньютона многочлена $f(z,w)-\xi_0$
совпадает с многоугольником Ньютона многочлена $f(z,w)-\xi$. Пункт
(Б) следует из следствия~\ref{newton:cor:1.3},
теоремы~\ref{newton:the:1.2} и следствия~\ref {newton:cor2}.
Пункт (В) следует из теоремы~\ref{newton:the:1.2} и следствия~\ref
{newton:cor2}.
\end{proof}

\section{Примеры} \label {sec:5}

\begin{figure}[h]
\begin{center}
\begin{minipage}[h]{0.49\linewidth}
\includegraphics[scale=0.6]{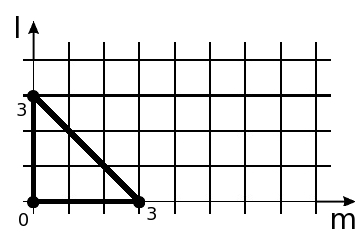}
\caption{Многоугольник Ньютона многочлена $f(z,w)=z^3+w^3-\xi$}
\label{newton:fig:ex2}
\end{minipage}
\hfill
\begin{minipage}[h]{0.49\linewidth}
\includegraphics[scale=0.6]{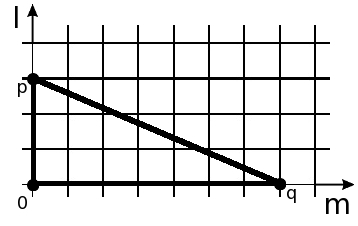}
\caption{Многоугольник Ньютона многочлена $f(z,w)=z^p+w^q-\xi$}
\label{newton:fig:ex3}
\end{minipage}
\end{center}
\end{figure}

\begin{example}
Пусть $f(z,w)=z^3+w^3$, тогда особое значение одно и равно $0$. Для
любого $\xi\in\CC\setminus\{0\}$ неособый слой
$\SSigma_\xi\approx\mathbb{T}^2\setminus\left\{p_{\xi,1},p_{\xi,2},p_{\xi,3}\right\}$
--- тор без трех бесконечно удаленных точек. В каждой бесконечно удаленной точке $p_{\xi,i}$, $i=1,2,3$,
векторное поле $\sgrad f$ имеет устранимую особенность.
Действительно, критическая точка равна $(0,0)$, особое значение равно
$0$, множество неособых значений совпадает с $\CC\setminus\{0\}$.
Многоугольник Ньютона для многочлена $f(z,w)-\xi=z^3+w^3-\xi$,
$\xi\in\CC\setminus\{0\}$, является треугольником с вершинами в
точках $A_1(3,0)$, $A_2(0,3)$, $A_3(0,0)$, см.
рис.~\ref{newton:fig:ex2}. Проверим, что многочлен $f(z,w)-\xi$,
$\xi\in\CC\setminus\{0\}$, является невырожденным относительно своего
многоугольника Ньютона. Обозначим сторону $A_2A_3$ многоугольника
Ньютона через $\Gamma_1$, сторону $A_1A_3$ через $\Gamma_2$, сторону
$A_1A_2$ через $\Gamma_3$. Многочлен, отвечающий граням $\Gamma_1$ и
$\Gamma_2$, равен $P_{\Gamma_1}(y)=P_{\Gamma_2}(y)=y^3-\xi$ и не
имеет кратных корней. Многочлен, отвечающий стороне $\Gamma_3$, равен
$P_{\Gamma_3}(y)=y^3+1$ и не имеет кратных корней. Отсюда следует,
что условия теоремы~\ref{newton:intrth2} выполнены. Количество
целочисленных точек строго внутри многоугольника Ньютона равно
$n_g=1$, количество целочисленных точек на стороне $\Gamma_3$ равно
2, поэтому $n_\mu=3$. По теореме~\ref{newton:intrth2} слой
$\SSigma_\xi$ имеет требуемые свойства.
\end{example}

\begin{example}
Пусть $f(z,w)=z^p+w^q$, $p,q\in\NN$, тогда особое значение одно и
равно $0$. Для любого $\xi\in\CC\setminus\{0\}$ неособый слой
$\SSigma_\xi$ гомеоморфен сфере с
$\left((p-1)(q-1)-(\gcd(p,q)-1)\right)/2$ ручками и без $\gcd(p,q)$
бесконечно удаленных точек. В каждой бесконечно удаленной точке
$p_{\xi,i}$, $i=1,\dots,\gcd(p,q)$, векторное поле $\sgrad f$ имеет
особенность полюс порядка $\frac{(p-1)(q-1)-1}{\gcd(p,q)}-1$.
Действительно, критическая точка равна $(0,0)$, особое значение равно
$0$, множество неособых значений совпадает с $\CC\setminus\{0\}$.
Многоугольник Ньютона для многочлена $f(z,w)-\xi=z^p+w^q-\xi$,
$\xi\in\CC\setminus\{0\}$, является треугольником с вершинами в
точках $A_1(q,0)$, $A_2(0,p)$, $A_3(0,0)$, см.
рис~\ref{newton:fig:ex3}. Проверим, что многочлен $f(z,w)-\xi$,
$\xi\in\CC\setminus\{0\}$, является невырожденным относительно своего
многоугольника Ньютона. Обозначим сторону $A_2A_3$ многоугольника
Ньютона через $\Gamma_1$, сторону $A_1A_3$ через $\Gamma_2$, сторону
$A_1A_2$ через $\Gamma_3$. Многочлен, отвечающий стороне $\Gamma_1$,
равен $P_{\Gamma_1}(y)=y^p-\xi$ и не имеет кратных корней. Многочлен,
отвечающий стороне $\Gamma_2$, равен $P_{\Gamma_2}(y)=y^q-\xi$ и не
имеет кратных корней. Многочлен, отвечающий стороне $\Gamma_3$, равен
$P_{\Gamma_3}(y)=y^{\gcd(p,q)}+1$ и не имеет кратных корней. Отсюда
следует, что условия теоремы~\ref{newton:intrth2} выполнены.
Количество целочисленных точек строго внутри многоугольника Ньютона
равно $n_g=((p-1)(q-1)-(\gcd(p,q)-1))/2$, количество целочисленных
точек на стороне $\Gamma_3$ равно $\gcd(p,q)-1$, поэтому
$n_\mu=\gcd(p,q)$. По теореме~\ref{newton:intrth2} слой $\SSigma_\xi$
имеет требуемые свойства.
\end{example}

\begin{figure}[h]
\begin{center}
\begin{minipage}[h]{0.49\linewidth}
\includegraphics[scale=0.6]{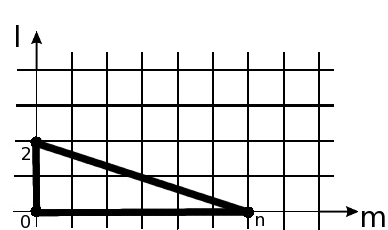}
\caption{Многоугольник Ньютона многочлена $f(z,w)=z^2+P_n(w)-\xi$
при $\xi\ne a_0$} \label{newton:fig:ex4}
\end{minipage}
\hfill
\begin{minipage}[h]{0.49\linewidth}
\includegraphics[scale=0.6]{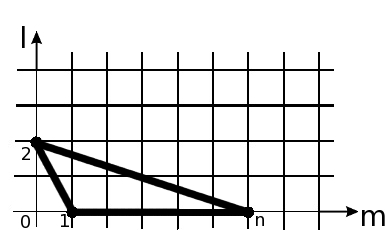}
\caption{Многоугольник Ньютона многочлена $f(z,w)=z^2+P_n(w)-\xi$
при $\xi=a_0$} \label{newton:fig:ex4:2}
\end{minipage}
\end{center}
\end{figure}

\begin{example} Пусть $f(z,w)=z^2+P_n(w)$, где $P_n(w)$ --- многочлен
одной переменной степени $n$, $P_n(w)=\sum\limits_{k=0}^{n}a_kw^k$,
$a_0\dots,a_n\in\CC$, $a_n\ne0$, тогда особые значения равны
$\xi_i=P_n(w^{0}_i)$, где $w^{0}_i$ --- корень уравнения $P'_n(w)=0$,
$i=1,\dots,n-1$. Для любого $\xi\in\CC\setminus\{\xi_i\}_{i=1}^{n-1}$
неособый слой $\SSigma_\xi$ гомеоморфен сфере с
$\left[(n-1)/2\right]$ ручками и без $\left(3+(-1)^n\right)/2$
бесконечно удаленных точек. В каждой бесконечно удаленной точке
$p_{\xi,i}$, $i=1,\dots,\left(3+(-1)^n\right)/2$, векторное поле
$\sgrad f$ имеет особенность полюс порядка $\frac{n-2}{\gcd(n,2)}-1$.
Действительно, критическая точка равна $(0,w^{0}_i)$, где $w^{0}_i$
--- корень уравнения $P'_n(w)=0$, $i=1,\dots,n-1$, отсюда особые
значения равны $\xi_i=f(0,w^{0}_i)=P_n(w^{0}_i)$, множество неособых
значений совпадает с $\CC\setminus\{\xi_i\}_{i=1}^{n-1}$.
Многоугольник Ньютона для многочлена $f(z,w)-\xi=z^2+P_n(w)-\xi$,
$\xi\in\CC\setminus\{\xi_i\}_{i=1}^{n-1}$, является треугольником с
вершинами в точках $A_1(n,0)$, $A_2(0,2)$, $A_3(0,0)$, если
$a_0-\xi\ne0$, см. рис~\ref{newton:fig:ex4}, и треугольником с
вершинами в точках $A_1(n,0)$, $A_2(0,2)$, $A_3(1,0)$, если
$a_0-\xi=0$, причем в этом случае $a_1\ne0$, так как $\xi$ ---
неособое значение, см. рис.~\ref{newton:fig:ex4:2}. Проверим, что
многочлен $f(z,w)-\xi$, $\xi\in\CC\setminus\{\xi_i\}_{i=1}^{n-1}$,
является невырожденным относительно своего многоугольника Ньютона.
Обозначим сторону $A_2A_3$ многоугольника Ньютона через $\Gamma_1$,
сторону $A_1A_3$ через $\Gamma_2$, сторону $A_1A_2$ через $\Gamma_3$.
Многочлен, отвечающий стороне $\Gamma_1$, равен либо
$P_{\Gamma_1}(y)=y^2+a_0-\xi$, если $a_0-\xi\ne0$, либо
$P_{\Gamma_1}(y)=y^2+a_1$, если $a_0-\xi=0$. В силу того, что $\xi$
--- неособое значение, $a_0-\xi=0$ и $a_1=0$ не могут выполняться
одновременно, поэтому $P_{\Gamma_1}$ не имеет кратных корней.
Многочлен, отвечающий стороне $\Gamma_2$, равен
$P_{\Gamma_2}(y)=P_n(y)-\xi$ и не имеет кратных корней. Многочлен,
отвечающий стороне $\Gamma_3$, равен $P_{\Gamma_3}(y)=y^2+a_n$ и не
имеет кратных корней. Отсюда следует, что условия
теоремы~\ref{newton:intrth2} выполнены. Количество целочисленных
точек строго внутри многоугольника Ньютона равно
$n_g=\left[(n-1)/2\right]$, количество целочисленных точек на стороне
$\Gamma_3$ равно $\left(1+(-1)^n\right)/2$, поэтому
$n_\mu=\left(3+(-1)^n\right)/2$. По теореме~\ref{newton:intrth2} слой
$\SSigma_\xi$ имеет требуемые свойства. \end{example}

\end{document}